\documentclass{amsart}
\usepackage{amssymb}
\usepackage{amsfonts,amsmath,amssymb,mathrsfs,verbatim}
\usepackage[dvips]{graphicx}
\usepackage{graphics}
\usepackage{psfrag}
\usepackage{hyperref}

\begin{document}
\newtheorem{thm1}{Theorem}[section]
\newtheorem{lem1}[thm1]{Lemma}
\newtheorem{rem1}[thm1]{Remark}
\newtheorem{def1}[thm1]{Definition}
\newtheorem{cor1}[thm1]{Corollary}
\newtheorem{defn1}[thm1]{Definition}
\newtheorem{prop1}[thm1]{Proposition}
\newtheorem{ex1}[thm1]{Example}
\newtheorem{alg1}[thm1]{Algorithm}


\title{Generalized robust toric ideals}
\author{Christos Tatakis}
\email{chtataki@cc.uoi.gr, chtatakis@gmail.com}
\address{Department of Mathematics, University of the Aegean, Samos 83200, Greece}

\begin{abstract}
\par
An ideal $I$ is robust if its universal Gr\"obner basis is a minimal generating set for
this ideal. In this paper, we generalize the meaning of robust ideals. An ideal is
defined as generalized robust if its universal Gr\"obner basis is
equal to its universal Markov basis. This article consists of two parts. In the first
one, we study the generalized robustness on toric ideals of a graph $G$. We prove that a
toric graph ideal is generalized robust if and only if its universal Markov basis is
equal to the Graver basis
of the ideal. Furthermore, we give a graph theoretical characterization of generalized
robust graph ideals, which is based on terms of graph theoretical properties of the
circuits of the graph $G$. In the second part, we go on to describe
the general case of toric ideals, in which we
prove that a robust toric ideal has a unique minimal system of generators, or in other words,
all of its minimal generators are indispensable.
\end{abstract}

\maketitle

\section{Introduction}

Let $A=\{\textbf{a}_1,\ldots,\textbf{a}_m\}\subseteq \mathbb{N}^n$
be a vector configuration in $\mathbb{Q}^n$ and
$\mathbb{N}A:=\{l_1\textbf{a}_1+\cdots+l_m\textbf{a}_m \ | \ l_i \in
\mathbb{N}\}$ the corresponding affine semigroup, where $\mathbb{N}A$
is pointed, that is if $x\in \mathbb{N}A$ and $-x\in \mathbb{N}A$ then $x={\bf{0}}$.
We grade the
polynomial ring $\mathbb{K}[x_1,\ldots,x_m]$ over an arbitrary field
$\mathbb{K}$ by the semigroup $\mathbb{N}A$ setting
$\deg_{A}(x_i)=\textbf{a}_i$ for $i=1,\ldots,m$. For
$\textbf{u}=(u_1,\ldots,u_m) \in \mathbb{N}^m$, we define the
$A$-degree of the monomial
$\textbf{x}^{\textbf{u}}:=x_1^{u_1} \cdots x_m^{u_m}$ to be \[
\deg_{A}(\textbf{x}^{\textbf{u}}):=u_1\textbf{a}_1+\cdots+u_m\textbf{a}_m
\in \mathbb{N}A.\]  The toric ideal $I_{A}$ associated to
$A$ is the prime ideal generated by all the binomials
$\textbf{x}^{\textbf{u}}- \textbf{x}^{\textbf{v}}$ such that
$\deg_{A}(\textbf{x}^{\textbf{u}})=\deg_{A}(\textbf{x}^{\textbf{v}})$,
see \cite{ST}.

Toric ideals consist a special class of ideals in a polynomial ring.
They define toric varieties, a large class of algebraic
varieties, that play an important role to the development
of mathematics the last years. Their study starts with Hochster in \cite{HO}
and spreads through a series of lectures by Fulton, see \cite{FU, FU2}.
 As far as the applicability of toric ideals is concerned,
 it has to be mentioned that toric ideals are related to recent advances in polyhedral geometry,
 toric geometry, algebraic geometry, algebraic statistic, integer programming, graph theory,
 computation algebra e.t.c., where they are applied in a natural way,
 see for example \cite{DST,ES,MS,ST}.

There are several sets for a toric ideal, which include
crucial information about it, such as the Graver basis,
the universal Markov basis, the universal Gr\"obner basis and
the set of the circuits.
An irreducible binomial $\textbf{x}^{\textbf{u}}-
\textbf{x}^{\textbf{v}}$ in $I_A$ is called primitive if
there is no other binomial
 $\textbf{x}^{\textbf{w}}- \textbf{x}^{\textbf{z}}$ in $I_A$,
such that $\textbf{x}^{\textbf{w}}$ divides $
\textbf{x}^{\textbf{u}}$ and $\textbf{x}^{\textbf{z}}$ divides $
\textbf{x}^{\textbf{v}}$. The set of primitive binomials forms the
Graver basis of $I_A$ and is denoted by $Gr_A$. As it is known
by a theorem of Diaconis and Sturmfels, every minimal generating set of $I_A$
corresponds to a minimal Markov basis of $A$, which is denoted
 by $M_A$, see \cite[Theorem 3.1]{DST}. The universal Markov
 basis of $A$ is denoted by $\mathcal{M}_A$ and is defined as
 the union of all minimal Markov bases of $A$, see \cite[Definition 3.1.]{HS}.
 The universal Gr\"{o}bner basis of an ideal $I_A$, which is denoted by $\mathcal{U}_A$,
  is a finite subset of $I_A$ and it is a Gr\"{o}bner basis for the ideal
 with respect to all admissible term orders, see \cite{ST}. The support of
 a monomial ${\bf{x^u}}$ of $\mathbb{K}[x_1,\ldots,x_m]$
is $\textrm{supp}({\bf{x^u}}):=\{i\ | \ x_{i}\ \mbox{ divides}\ {\bf{x^u}}\}$ and
the support of a binomial $B={\bf{x^u}}-{\bf{x^v}}$ is
$\textrm{supp}(B):=\textrm{supp}({\bf{x^u}})\cup \textrm{supp}({\bf{x^v}})$.
An irreducible binomial is called circuit if it has minimal support.
The set of the circuits of a toric ideal $I_A$ is denoted by
$\mathcal{ C}_A$. The relation between the above sets was studied by B.~Sturmfels
in \cite{ST}:

\begin{prop1}\label{CUG}\cite[Proposition 4.11]{ST} For any toric ideal $I_A$ it holds:
$$\mathcal{ C}_A\subseteq \mathcal{ U}_A \subseteq  Gr_A.$$
\end{prop1}

An ideal $I$ is called robust if its universal Gr\"obner
basis is equal with  a Markov basis of the ideal. Robustness
is a property of ideals that has not been fully described. More specifically, it has been
described for toric ideals which are generated by quadratics. Toric ideals
which are generated by quadratics were studied by Ohsugi and Hibi in \cite{OH}, while
the robustness for this class of ideals is described in the article of Boocher and Robeva,
see \cite{BOO}. The importance of robustness stems from the interest
in the study of ideals which are minimally generated by a Gr\"obner basis for an arbitrary term
order, see \cite{CHT}. Moreover, the study of robustness is important, due to the fact
that several areas of mathematics are keen on the research of the
Markov basis, the universal Gr\"obner basis and the Graver basis of an ideal.
This problem has also been researched in the case of toric ideals
arising from a graph $G$, as studied by Boocher et al in \cite{BOO2}. In their work
the authors proved that any robust toric ideal of a graph $G$ is also minimally
generated by its Graver basis, \cite[Theorem 3.2.]{BOO2}. In addition, they completely characterize
all graphs which give rise to robust ideals, see \cite[Theorem 4.8.]{BOO2}.\par
The present article generalizes the meaning of robust ideals. A robust ideal is called
generalized robust if its universal Gr\"obner basis is equal with its universal Markov
basis. This manuscript is divided into two parts.\par
In the first part, we study the generalized robustness on toric ideals of a graph $G$.
The results of this part are inspired and guided by the work of
\cite{BOO2} in order to give theorems that fully characterize the
generalized robust toric ideals of graphs. The papers \cite{TT1},
\cite{TT2} and \cite{VI} describe
the Markov basis, the Graver basis, the universal Gr\"obner basis and the
set of the circuits for a toric ideal arising from a graph.
In section 2, we analyze all these notions more explicitly. Applying
this knowledge on the work of Boocher et al (see \cite{BOO2}), we are allowed to
provide the study of the generalized robustness of graphs,
 with theorems of the same structure as theirs. In section 3, we first prove that a
toric graph ideal is generalized robust if and only if its universal Markov basis is
equal to the Graver basis
of the ideal, see Theorem \ref{M=Gr}. Moreover, the
relation between robust graph ideals and
 generalized robust graph ideals is studied. In the next section, we go on
to give a graph theoretical characterization of generalized
robust graph ideals, which is based on terms of graph theoretical properties of the
circuits of the graph $G$, see Theorem \ref{circ-gen}.\par
 In the second part of this manuscript, we study the robustness property in the general case of toric ideals. More
 especially, we study
 the indispensable binomials which exist in a robust toric ideal. A binomial $B\in I_A$
 is called indispensable if there exists a non zero constant multiple of it
 in every minimal system of binomial generators of $I_A$. A recent problem arising
 from algebraic statistics is to find classes of toric ideals which have a unique
 minimal system of generators, see \cite{Tak},\cite{THO3}. In order to study this
 problem, Ohsugi and Hibi introduced in \cite{OH2} the notion of indispensable
binomials. In section 5, we
prove that a robust toric ideal has a unique minimal system of generators, or in other words,
all of its minimal generators are indispensable, see Theorem \ref{main}. Finally,
we conclude that the robustness property for a toric ideal, implies
the generalized robustness property for it, see Corollary \ref{eis}. In conclusion,
we present a family of toric ideals, for which $\mathcal{C}_A=Gr_A$, see Remark \ref{remm}.

\section{Elements of the toric ideals of graphs}

In the next chapters,  $G$ is a connected, undirected, finite, simple graph on the vertex set
$V(G)=\{v_{1},\ldots,v_{n}\}$.  Let $E(G)=\{e_{1},\ldots,e_{m}\}$ be
the set of edges of $G$ and $\mathbb{K}[e_{1},\ldots,e_{m}]$
 the polynomial ring in the $m$ variables $e_{1},\ldots,e_{m}$ over a field $\mathbb{K}$.  We
will associate each edge $e=\{v_{i},v_{j}\}\in E(G)$  with the element
$a_{e}=v_{i}+v_{j}$ in the free abelian group $ \mathbb{Z}^n $, with basis the set of the vertices
of $G$, where $v_{i}=(0,\ldots,0,1,0,\ldots,0)$ be the vector with 1 in the $i-$th coordinate of $v_{i}$. With $I_{G}$ we denote
 the toric ideal $I_{A_{G}}$ in
$\mathbb{K}[e_{1},\ldots,e_{m}]$, where  $A_{G}=\{a_{e}\ | \ e\in E(G)\}\subset \mathbb{Z}^n $.

\par
A walk of length $q$ connecting $v_{i_1}\in V(G)$ and $v_{i_{q+1}}\in V(G)$ is a finite sequence of
the form $w=(e_{i_1}=\{v_{i_1},v_{i_2}\},\{v_{i_2},v_{i_3}\},\ldots,e_{i_q}=\{v_{i_q},v_{i_{q+1}}\})$
with each $e_{i_j}=\{v_{i_j},v_{i_{j+1}}\}\in E(G)$. We call a walk $w'=(e_{j_{1}},\dots,e_{j_{t}})$
a subwalk of $w$ if $e_{j_1}\cdots e_{j_t}|\ e_{i_1}\cdots e_{i_q}.$ An
even (respectively odd) walk is a walk of even (respectively odd) length.
A walk as $w$ is called closed if $v_{i_{q+1}}=v_{i_1}$. A cycle
is a closed walk with $v_{i_k}\neq v_{i_j},$ for every $ 1\leq k < j \leq q$. Depending
on the property of the walk that we want to emphasize, we may
denote a walk $w$ either by a sequence of vertices and edges $(v_{i_1},
e_{i_1}, v_{i_2}, \ldots ,v_{i_q}, e_{i_q}, v_{i_{q+1}})$ or exclusively with
vertices $(v_{i_1},v_{i_2},v_{i_3},\ldots,v_{i_{q+1}})$ or only with
edges $(e_{i_1},\ldots ,e_{i_q})$. Note that, although the graph
$G$ has no multiple edges, since it is simple, the
same edge $e$ may appear more than once in a walk. In this case, $e$ is
called multiple edge of the walk $w$.
Given an even closed walk $w =(e_{i_1}, e_{i_2},\dots,
e_{i_{2q}})$ of the graph $G$, we denote by $B_w$ the binomial
$$B_w=\prod _{k=1}^{q} e_{i_{2k-1}}-\prod _{k=1}^{q} e_{i_{2k}}$$
belonging to the toric ideal $I_G$. Actually, the toric ideal $I_G$
is generated by binomials of the above form, see \cite{VI}.

For convenience, by $\bf{w}$
we denote the subgraph of $G$, whose vertices and edges
 are the vertices and the edges of the walk $w$. Note that $\bf{w}$
is a connected subgraph of $G$.
A cut edge (respectively a cut vertex) is an edge (respectively a vertex) of
the graph, whose removal increases the number of connected
components of the remaining subgraph.  A graph is called biconnected
if it is connected and does not contain a cut
vertex. A block is a maximal biconnected subgraph of a given
graph $G$.

A walk $w$ of a graph is primitive if and only if the corresponding binomial
$B_w$ is primitive. Every even primitive walk $w=(e_{i_1},\ldots,e_{i_{2k}})$
partitions the set of its edges in two sets $w^+= \{e_{i_j}|j \
 \textrm{odd}\}$ and $w^-=\{e_{i_j}|j \ \textrm{even}\}$, otherwise the
binomial $B_w$ is not irreducible. The edges of the set $w^+$ are called
odd edges of the walk and those of
$w^-$ even. A sink of a block $\mathcal{B}$ is a
common vertex of two odd or two even edges of the walk $w$ which
belong to the block $\mathcal{B}$. Finally,
we call strongly primitive walk, a primitive walk which has not
two sinks with distance one in any cyclic block of the walk, or
equivalently has not two adjacent cut vertices in any cyclic block of $w$.

In \cite{TT1} a complete
characterization of the Graver basis of the corresponding toric ideal $I_G$ was given,
see \cite[Theorem 3.2]{TT1}. The next corollary, given by the same authors,
describes the structure of the underlying graph of a primitive walk.

\begin{cor1}\cite[Corollary 3.3]{TT1} \label{primitive-graph}
Let $G$ be a graph and $W$ a connected subgraph of $G$. The subgraph
$W$ is the graph  ${\bf w}$ of a primitive walk $w$ if and only if
\begin{enumerate}
  \item  $W$ is an even cycle or
  \item  $W$ is not biconnected and
\begin{enumerate}
  \item every block of $W$ is a cycle or a cut edge and
  \item every cut vertex of $W$ belongs to exactly
  two blocks and separates the graph in two parts, the total number of
edges of the cyclic blocks in each part is odd.
\end{enumerate}
\end{enumerate}
\end{cor1}

Afterwards, we recall from \cite{TT1}, a lot of graph theoretical notions
 in order to describe the universal Markov basis of a toric ideal of a graph $G$.
 We say that a binomial is a minimal binomial,
if it belongs to at least one minimal system of generators of $I_G$, i.e.
at least one Markov basis of $I_G$.

For a given subgraph $F$  of $G$, an edge $f$
of the graph $G$ is called chord of the subgraph $F$, if the vertices of the edge
$f$ belong to $V(F)$ and $f\notin E(F)$. A chord $e=\{v_k,v_l\}$ is called
bridge of a primitive walk $w$ if there exist two different blocks
${\mathcal B}_1,{\mathcal B}_2$ of $\bf{w}$ such that $v_k\in {\mathcal B}_1$
and $v_l\in {\mathcal B}_2$. Let $w$ be an even closed walk
$(\{v_{1},v_{2}\},\{v_{2},v_{3}\},\ldots,\{v_{2q},v_{1}\})$ and
$f=\{v_{i},v_{j}\}$ a chord of $w$. Then, $f$ breaks $w$ into two
walks:
$$w_{1}=(e_{1},\ldots,e_{i-1}, f, e_{j},\ldots,e_{2q})$$ and
$$w_{2}=(e_{i},\ldots,e_{j-1},f),$$ where $e_{s}=\{v_{s},v_{s+1}\},\ 1\leq s\leq 2q$ and $e_{2q}=\{v_{2q},v_{1}\}.$
The two walks are both even or both odd. A chord is called even (respectively odd) if it is not
a bridge and if it breaks the walk into two even walks (respectively odd).

Let
$w=(\{v_{i_{1}},v_{i_{2}}\}, \{v_{i_{2}},v_{i_{3}}\},\cdots ,
\{v_{i_{2q}},v_{i_{1}}\})$ be a primitive walk. Let
$f=\{v_{i_{s}},v_{i_{j}}\}$ and $f'=\{v_{i_{s'}},v_{i_{j'}}\}$ be two
odd chords (that is they are not bridges and the numbers
$j-s,j'-s'$ are even) with $1\leq
s<j\leq 2q$ and $1\leq s'<j'\leq 2q$. We say that $f$ and $f'$
cross effectively in $w$ if $s'-s$ is odd (then necessarily $j-s',
j'-j, j'-s$ are odd) and  either $s<s'<j<j'$ or $s'<s<j'<j$.
We call $F_4$ of the walk $w$, a cycle $(e, f, e', f')$ of
length four which consists of two edges $e,e'$ of the walk $w$
both odd or both even, and two odd chords  $f$,$f'$ which
cross  effectively in $w$. An $F_4$, $(e_1, f_1, e_2, f_2)$ separates the vertices of
${\bf w}$ into two parts $V({\bf w}_1), V({\bf w}_2)$,
since both edges $e_1, e_2$  of the $F_4$ belong to the same
block of $w=(w_1, e_1, w_2, e_2)$. We say that an odd chord $f$ of
a primitive walk $w=(w_1, e_1, w_2, e_2)$ crosses an $F_4$, $(e_1, f_1, e_2, f_2)$,
if  one of the vertices of $f$ belongs to $V({\bf w}_1)$, the
other in $V({\bf w}_2)$ and $f$ is different from  $f_1$ and $f_2$.

The next theorem by Reyes et al, from
which we know the elements of the universal Markov basis of the ideal $I_G$,
gives a necessary and sufficient
characterization of minimal binomials of a toric ideal of a graph $G$.

\begin{thm1}\cite[Theorem 4.13]{TT1} \label{minimal}
Let $w$ be an even closed walk. $B_{w}$ is a minimal binomial if
and only if \begin{enumerate}
  \item[(M1)] all the chords of $w$ are odd,
  \item[(M2)] there are not two odd chords of $w$ which cross effectively except if they form an $F_4$,
  \item[(M3)] no odd chord crosses an $F_4$ of the walk $w$,
  \item[(M4)] $w$ is a strongly primitive.
\end{enumerate}
\end{thm1}

\section{Generalized robust toric ideals of graphs}

In this section we study the generalized robust toric ideals.

\begin{def1}An ideal $I$ is called generalized robust if its universal Gr\"obner
basis is equal with its universal Markov basis.
\end{def1}

Undoubtedly, it is a hard problem to characterize the generalized robustness, owing to the
fact that only for a few classes of toric ideals we know their universal
Gr\"obner basis and the universal Markov basis. In general,
characterizing and computing these sets, is a difficult and computationally demanding
problem.
Lawrence ideals
 provide a large class of generalized robust toric ideals, since
 it is known by Sturmfels that in a Lawrence ideal any minimal
 generating set coincides with the universal
 Gr\"obner basis and the Graver basis, see \cite[Theorem 7.1]{ST}.
 Moreover, not only the robustness but also the generalized
 robustness is not a property which describes completely
  the Lawrence ideals, see \cite[Example 3.4]{BOO}.

In order to describe the universal
Gr\"obner basis for the case of toric ideals of graphs, we give the notions of pure blocks
 and of the mixed walks of a graph $G$, see \cite{TT2}.
 A cyclic block $\mathcal{B}$ of a primitive walk $w$ is called pure if all the edges
of the block $\mathcal{B}$ belong either to $\textbf{w}^+$ or to $\textbf{w}^-$. A primitive
walk $w$ is called mixed if none of the cyclic blocks of $w$ is pure. The next theorem
describes completely the elements of the universal Gr\"obner
basis of a toric ideal of a graph $G$.
\begin{thm1}\cite[Theorem 3.4]{TT2} \label{UGB} Let $w$ be a primitive walk.
$B_w$ belongs to the universal Gr\"{o}bner basis of $I_G$ if and only if $w$ is mixed.
\end{thm1}

Based on the above theorem, in combination with the knowledge of the universal Markov basis
for a toric ideal of a graph $G$, we are allowed
to research in depth the generalized robust graph ideals.

In the special case of toric ideals of graphs, a useful property for every minimal generator of the ideal,
is that it
belongs to its universal Gr\"obner basis, as we can see in the next proposition.

\begin{prop1} \label{M<U}Let $G$ be a graph and $I_G$ its corresponding toric ideal.
Then $$\mathcal{M}_G\subseteq\mathcal{U}_G.$$
\end{prop1}
\textbf{Proof.} Let $G$ be a graph and $I_G$ its corresponding
toric ideal. Let $B_w$ be an element of the universal Markov basis of $I_G$,
 which means by definition that $B_w$ is a minimal generator of
 the ideal. We will prove that the binomial $B_w$ belongs to the universal
 Gr\"obner basis of $I_G$. \par We assume that $B_w$ does not
  belong to the $\mathcal{U}_G$. By Theorem \ref{UGB} the walk $w$ is not mixed.
 Therefore, the walk $w$ has at least one pure cyclic block and let it be
 $\mathcal{B}=(e_1,\ldots,e_n)$. Thus, all the edges $e_i$ of the block $\mathcal{B}$
 are either even or odd. Let $e_1=(u_1,u_2),e_2=(u_2,u_3)$ and $e_3=(u_3,u_4)$ be
 three consecutive edges of the block $\mathcal{B}$. We know that the block
 $\mathcal{B}$ has at least three edges, since
 $\mathcal{B}$ is a cycle of the graph $G$. The vertices $u_2$ and $u_3$ are both common vertices of either two odd
 or two even edges of $\mathcal{B}$. In consequence, the vertices $u_2$ and $u_3$ are sinks of the walk $w$
 with distance one. We remark that $w$ is primitive, since $B_w$
 is a minimal generator of the ideal. It follows that $w$ is not strongly primitive,
a contradiction arises due to the minimality of $B_w$ and Theorem \ref{minimal}.
\hfill $\square$

We mention that the above argument is not true in the general case of toric ideals.
In \cite{THO1} the authors provide a counterexample for this claim,
\cite[Example 1.8.]{THO1}.

We are now ready to prove our main result in this section, in which we describe the generalized
robustness for graph ideals. In \cite{BOO2} the authors
proved that a graph ideal is robust
if and only if the Graver basis of the ideal is equal to a Markov
basis of it. Next, we are proving the corresponding theorem for
a generalized robust ideal which is stated on its Graver basis and its universal
Markov basis of the ideal.
\begin{thm1}\label{M=Gr}Let $G$ be a graph and $I_G$ its corresponding ideal. The ideal $I_G$ is
generalized robust if and only if $\mathcal{M}_G=Gr_G$.
\end{thm1}
\textbf{Proof.} Let $G$ be a graph and $I_G$ its corresponding toric ideal.
From Proposition \ref{CUG} and Proposition \ref{M<U} we have
that $$\mathcal{M}_G\subseteq \mathcal{U}_G\subseteq Gr_G.$$
If $\mathcal{M}_G=Gr_G$, then the result follows.
\par
\begin{figure}[here]
\begin{center}
\psfrag{A}{$w_{1}$}\psfrag{B}{$w_{2}$}\psfrag{C}{$w_{n}$}
\psfrag{D}{$e_{1}$}\psfrag{E}{$e_{2}$}\psfrag{F}{$e_{n}$}
\includegraphics{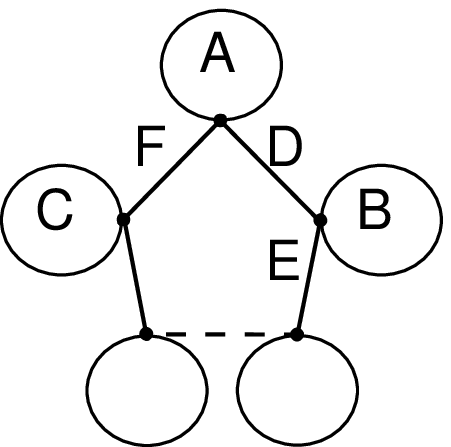}
\caption{The walk w}
\label{Figure 5}
\end{center}
\end{figure}

 Conversely, let $I_G$ be a generalized robust ideal,
which means that $\mathcal{M}_G=\mathcal{U}_G$. It is enough to prove
that $Gr_G\subseteq \mathcal{U}_G$. Let $w$ be a walk of the graph $G$
such that $B_w\in Gr_G$. We will prove that $B_w\in \mathcal{U}_G$. Suppose this is not true.
Therefore, the walk $w$ has at least one pure cyclic block, and let it
$\mathcal{B}=(e_1,\ldots,e_n)$. The walk $w$ has the form $w=(w_1,e_1,w_2,e_2,\ldots,e_{n-1},w_{n},e_n)$,
where $w_1,\ldots,w_n$ are odd subwalks of $w$, as in Figure \ref{Figure 5}. We remark
that the walks $w_1,\ldots,w_n$ are odd. The reason is that the walk $w$ is primitive, which means
that every cut vertex of $w$ separates the graph into two parts, the total number of
edges of the cyclic blocks in each part is odd, see Corollary \ref{primitive-graph}.

Furthermore, from Corollary \ref{primitive-graph} we know that every block of $w$ is
a cycle or a cut edge. Therefore, we can assume that each one of the walks $w_1$ and
$w_n$ has at least one odd cycle. Let them be $c_1$ and $c_n$ correspondingly.
We will prove the existence of a walk $q$
of the graph $G$, such that $B_q$ belongs to $\mathcal{U}_G$ but not
to $\mathcal{M}_G$.

We consider the walk
$$q=(c_1,p_1,e_1,e_2,\ldots,e_{n-1},p_2,c_n,-p_2,e_{n-1},\ldots,e_1,p_1)$$ which consists of the two
odd cycles $c_1,c_n$ and the path $p=(p_1,e_1,\ldots,e_{n-1},p_2)$ joining them, where
$p_1,p_2$ are paths of $w$ which join the cycles $c_1$ and $c_n$
with the edges $e_1$ and $e_{n-1}$ correspondingly. From Theorem
\ref{UGB} it follows that the binomial $B_q$ is an element of the $\mathcal{U}_G$.
 We remark that the edge $e_n$ of $w$ is a bridge of the walk $q$.
From Theorem \ref{minimal} it follows that the binomial $B_q$ is not minimal.
As a result, the binomial $B_q$ does not belong to the universal Markov basis
of the ideal $I_G$, a contradiction arises. \hfill $\square$

From the above theorem it follows that the robustness implies
 the generalized robustness for a toric graph ideal.

\begin{cor1}\label{rob-gen} Let $I_G$ be a robust ideal of a graph $G$. The ideal $I_G$ is generalized robust.
\end{cor1}

The converse of the above corollary is not true as the following
 example proves.
\begin{ex1}\label{ex1}{\rm
We consider the complete graph $G=K_4$ (see Figure \ref{complete}) on the four vertices and let $I_G$
be its corresponding toric ideal.

\begin{figure}[here]
\begin{center}
\psfrag{A}{$e_{1}$}\psfrag{B}{$e_{2}$}\psfrag{C}{$e_{3}$}\psfrag{D}{$e_{4}$}\psfrag{E}{$e_{6}$}\psfrag{F}{$e_{5}$}
\includegraphics{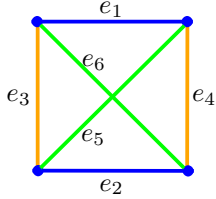}
\caption{The complete graph $K_4$}
\label{complete}
\end{center}
\end{figure}

The ideal has three minimal generators, which
are $$B_{w_1}=e_1e_2-e_5e_6, B_{w_2}=e_3e_4-e_5e_6\ \textrm{and}\ B_{w_3}=e_1e_2-e_3e_4.$$
Moreover, it has three minimal system of generators
$$M_1=<B_{w_1},B_{w_2}>, M_2=<B_{w_1},B_{w_3}>\ \textrm{and}\ M_3=<B_{w_2},B_{w_3}>,$$ and as a result its
universal Markov basis is $\mathcal{M}_G=<B_{w_1},B_{w_2},B_{w_3}>$.
Since the walks $w_1,w_2,w_3$ are even cycles and no other block exists
on the graph $G$, from Corollary \ref{primitive-graph}
it follows that there is no other primitive elements of $I_G$
and therefore $Gr_G=<B_{w_1},B_{w_2},B_{w_3}>$. Thus, the ideal $I_G$ is generalized
robust but not robust.
}
\end{ex1}

In order to check the converse statement of the last corollary, we will
use the following corollary, as it was presented in \cite{TT1}.

\begin{cor1}\cite[Corollary 4.15]{TT1}\label{uniqueminm} Let $G$ be a graph which has no cycles
of length four. The toric ideal $I_G$ has a unique system of binomial generators.
\end{cor1}

\begin{prop1} Let $I_G$ be a generalized robust ideal. If the graph $G$ has no cycles of
length four, then $I_G$ is robust.
\end{prop1}
\textbf{Proof.} Let $I_G$ be a generalized robust ideal, where $G$ has
no cycles of length four. By definition we know that $\mathcal{M}_G=\mathcal{U}_G$.
From Corollary \ref{uniqueminm} it follows that $I_G$ has a unique system of
minimal generators, which means that $M_G=\mathcal{M}_G$. Therefore, the ideal
$I_G$ is minimally generated by its Gr\"obner basis. It follows that $I_G$ is robust. \hfill$\square$

The converse of the above proposition is not true as the following example proves.
\begin{ex1}{\rm
Let $G$ be
the graph which is a chordless cycle of length four, see Figure \ref{square}.

\begin{figure}[here]
\begin{center}
\psfrag{A}{$e_{1}$}\psfrag{B}{$e_{2}$}\psfrag{C}{$e_{3}$}\psfrag{D}{$e_{4}$}
\includegraphics{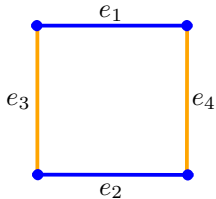}
\caption{A both robust and generalized robust graph}
\label{square}
\end{center}
\end{figure}

Then, the corresponding
toric ideal $I_G$ has one minimal generator the $B_w=<e_1e_2-e_3e_4>$.
Obviously, the ideal has a unique system of minimal generators, the $M_G=<B_w>$
and therefore $\mathcal{M}_G=<B_w>$. It is clear that the Graver basis of $I_G$
consists of exactly the binomial $B_w$. It follows that the ideal
is both robust and generalized robust, but it contains a cycle of length four.
}
\end{ex1}

Certainly, the uniqueness of the minimal system of generators of an ideal is
 a sufficient and a necessary condition for the toric ideal of a graph $G$,
 to be robust if it is generalized robust and conversely.

\section{Circuits and generalized robust graph ideals}

In this section, we will present a graph theoretical characterization
of a generalized robust toric ideal of a graph $G$, which is based on terms
of graph theoretical properties of the circuits of the graph $G$.

The
circuits of a graph $G$ were described in graph theoretical terms
with necessary and sufficient conditions by R.~Villarreal:
\begin{prop1}\cite[Proposition 4.2]{VI}\label{circuit}
Let $G$ be a finite connected graph. The binomial $B\in I_{G}$ is a circuit if and only if $B=B_{w}$ where
\begin{enumerate}
  \item[(C1)] $w$ is an even cycle or
  \item[(C2)] $w$ consists of two odd cycles intersecting in exactly one vertex or
  \item[(C3)] $w$ consists of two vertex disjoint odd cycles joined by a path.
\end{enumerate}
\end{prop1}

From \cite{OH} we also know the form of the primitive walks of a graph $G$.

\begin{lem1}\label{primitOH}\cite[Lemma 3.2]{OH} If $B_w$ is primitive, then $w$ has one of the
following forms:
\begin{enumerate}
  \item[(P1)] $w$ is an even cycle or
  \item[(P2)] $w$ consists of two odd cycles intersecting in exactly one vertex or
  \item[(P3)] $w=(c_1,w_1,c_2,w_2)$ where $c_1,c_2$ are odd vertex disjoint cycles and $w_1,w_2$
  are walks which combine a vertex $v_1$ of $c_1$ and a vertex $v_2$ of $c_2$.
\end{enumerate}
\end{lem1}

As we see in Corollary \ref{primitive-graph}, if a walk $w$ has one of the first two forms it is also primitive.
However, this is not true when the walk $w$ has the third form. As we saw, the Corollary \ref{primitive-graph}
describes completely the primitive graphs. It is clear
that in the case that the walk $w$ has either (P1) either (P2) form,
then it is always strongly primitive, since the corresponding primitive graph ${\bf{w}}$
does not contain two sinks.

The next proposition describes some properties of the primitive elements of a
generalized robust graph ideal which will be will be our main tool, to prove the main result in this section.

\begin{prop1}\label{propgen}The ideal $I_G$ is generalized robust if and only if all its primitive
elements satisfy the conditions M1 and M2 of the Theorem \ref{minimal}.
\end{prop1}
\textbf{Proof.} Let $I_G$ be a generalized robust ideal. From Theorem \ref{M=Gr}
all primitive elements are minimal generators and therefore they satisfy the conditions
$M1$ and $M2$. Conversely, we assume that
all the primitive elements of $I_G$ satisfy the conditions $M1$ and $M2$ of the
Theorem \ref{minimal}. In order to prove that the ideal $I_G$ is generalized robust,
from Theorem \ref{M=Gr} and Proposition \ref{M<U} we have to prove that $Gr_G\subseteq \mathcal{M}_G$.
Let $B_w$ be an element of $Gr_G$ such that it satisfies the conditions $M1$ and $M2$ of the
Theorem \ref{minimal}. We have to prove that the binomial $B_w$ is minimal.
By hypothesis, it remains to prove that $B_w$ satisfies the conditions $M3$ and $M4$ of the
Theorem \ref{minimal}, i.e. the walk $w$ has not an
odd chord which crosses an $F_4$ of $w$ and the walk $w$ is strongly primitive.

\begin{figure}[here]
\begin{center}
\psfrag{A}{$e_{2}$}\psfrag{B}{$e_{1}$}\psfrag{C}{$f_{1}$}\psfrag{D}{$f_{2}$}\psfrag{E}{$f$}\psfrag{F}{$v_{1}$}
\psfrag{H}{$v_{2}$}\psfrag{I}{$w_{1}$}\psfrag{J}{$w_{2}$}\psfrag{K}{$w'_{1}$}\psfrag{L}{$w''_{1}$}
\psfrag{M}{$w'_{2}$}\psfrag{N}{$w''_{2}$}
\includegraphics{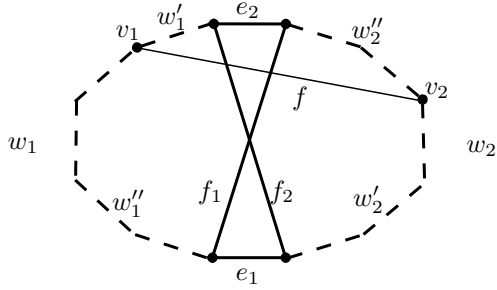}
\caption{An odd chord which crosses an $F_4$}
\label{f4}
\end{center}
\end{figure}

Firstly, we prove that $w$ has not an
odd chord which crosses an $F_4$ of $w$. Suppose not. So, there exists an odd chord $f=\{v_1, v_2\}$
that crosses the $F_4$, $(e_1, f_1, e_2, f_2)$ of the walk
$w=(w_1, e_1, w_2, e_2)$, see Figure \ref{f4}.

Then, $w$ can be written in the form
$(w_1',\{v_1\},w_1'', e_1, w_2',\{v_2\}, w_2'', e_2)$, where $w_1=(w_1',\{v_1\},w_1'')$
and $w_2=(w_2',\{v_2\}, w_2'')$. Since the chord $f$ is odd,
 by definition the walks $(f, w_2'',e_2, w_1')$ and $(f, w_1'',e_1, w_2')$
are both odd. In addition, since $(e_1, f_1, e_2, f_2)$ is an $F_4$,
the walks $w_1$ and $w_2$ are both odd. Therefore, $(w_1'', f_1, -w_2'',f)$
and $(w_1', f, -w_2',f_2)$ are both even. So, from the definition, $f$
is an even chord of $w'=(w_1, f_1, -w_2, f_2)$. Since the walk $w$ is primitive,
the walk $w'$ is also primitive
 and $w'$ has an even chord,
 a contradiction arises due to the fact that by hypothesis all primitive elements satisfy
the condition $M1$ of the Theorem \ref{minimal}.

It remains to prove that $w$ is strongly primitive. Suppose not.
Therefore the walk $w$ has the form $w=(c_1,w_1,c_2,w_2)$, where
$c_1,c_2$ are odd vertex disjoint cycles and $w_1,w_2$
  are walks which combine a vertex $v_1$ of $c_1$ and a vertex $v_2$ of $c_2$.
  Since the walk is not strongly primitive,
 there is a cyclic block $\mathcal{B}$ of the primitive walk $w$ in which there are
 two adjacent cut vertices $u_1$ and $u_2$ of $w$. Let the block be $\mathcal{B}=(q,e)$,
 where $e$ is the edge $\{u_1,u_2\}$ and $q$ be the path of the cyclic block $\mathcal{B}$
 which connects the vertices $u_1$ and $u_2$. Then $w$ can be written in the form
 $$w=(c_1,\{v_1\},w',\{u_1\},q,\{u_2\},w'',\{v_2\},c_2,\{v_2\},-w'',\{u_2,u_1\}=e,-w'),$$
 where $w_1=(\{v_1\},w',\{u_1\},q,\{u_2\},w'',\{v_2\})$ and $w_2=(\{v_2\},-w'',e,-w')$.
 Since the graph is connected, there is a path between any two vertices of the graph.
 Let $p_1$ be a path between the vertices $v_1$ and $u_1$ and let
 $p_2$ be a path between the vertices $v_2$ and $u_2$.
 Without loss of generality, we assume that the paths $p_1$ and $p_2$ are subwalks of the
 subwalks $w'$ and $w''$ of the walk $w$.
 We consider the walk $$w_c=(c_1,p_1,\{u_1\},q,\{u_2\},-p_2,c_2,p_2,\{u_2\},-q,\{u_1\},-p_1).$$
 The above walk is a circuit since it is in the form $(C3)$, where $c_1,c_2$ are the two
 vertex disjoint odd cycles and $p=(p_1,q,-p_2)$ the path which joins them.
 Therefore, the binomial $B_{w_c}$ is primitive and the edge $e=\{u_1,u_2\}$ is a bridge
 of $w_c$. By hypothesis all the primitive elements of $I_G$ satisfy the condition $M1$ of the
Theorem \ref{minimal}, which means that they have no bridges. A contradiction arises. \hfill$\square$

In \cite[Theorem 4.8]{BOO2} the authors proved the following lemma:

\begin{lem1}\label{texn2} Let $I_G$ be an ideal such that all the chords of all of its primitive elements are odd.
Then, there is no circuit of $G$ which shares exactly one edge (and no other vertices)
 with another circuit such that the shared edge is part of a cyclic block in both circuits.
\end{lem1}

We are ready to present the main result of this section. The only difference between it
and the corresponding result of Boocher et al in \cite{BOO2}, in the case of robust graph ideals,
 is that for an ideal $I_G$
to be generalized robust we allow to the
circuits of the graph the existence of two odd chords of the walk $w$ which form an $F_4$.
As a result of the previous proposition, the proof of the following
theorem respects completely the construction of the corresponding proof of
Theorem 4.8. in \cite{BOO2}. In this part of the proof
we refer to the corresponding proof of them.

\begin{thm1}\label{circ-gen} The ideal $I_G$ is generalized robust if and only if the following
conditions are satisfied.
\begin{enumerate}
  \item[(R1)] No circuit of $G$ has either an even chord or a bridge,
  \item[(R2)] No circuit of $G$ contains two odd chords which cross effectively, except if
  they form an $F_4$,
  \item[(R3)] No circuit of $G$ shares exactly one edge (and no other vertices) with
  another circuit such that the shared edge is a part of a cyclic block in both circuits.
\end{enumerate}
\end{thm1}
\textbf{Proof.} By Proposition \ref{propgen} it is equivalently to prove that all
primitive elements of $I_G$ satisfy the conditions $M1$ and $M2$ of the Theorem \ref{minimal}
if and only if the circuits of the graph $G$ satisfy the conditions $R1$ through $R3$.\\
For the forward direction we assume that all primitive elements of $I_G$ satisfy
the conditions $M1,M2$ of the Theorem \ref{minimal}. Since all the circuits of $G$ are
also primitive elements, the conditions $R1,R2$ are followed. By Lemma \ref{texn2} it
follows the condition $R3$.

For the other direction, we assume that every circuit of $G$ satisfies $R1$ through $R3$. We will
prove that every primitive element of $I_G$ satisfies the conditions $M1$ and $M2$
of the Theorem \ref{minimal}. Suppose not. Then there exists a primitive walk $w$ of $I_G$
such that it has either an even chord or a bridge or it has two odd chords which
cross effectively and they do not form an $F_4$. From Lemma \ref{primitOH} the walk
$w$ is either an even cycle or two odd cycles intersecting in exactly one vertex or
 $w=(c_1,w_1,c_2,w_2)$ where $c_1,c_2$ are odd vertex disjoint cycles and $w_1,w_2$
  are walks which combine a vertex $v_1$ of $c_1$ and a vertex $v_2$ of $c_2$. If
  $w$ has one of the first two forms, then the walk $w$ is also a circuit,
   contradicted the hypothesis. Therefore, the walk $w$ is of the form
   $w=(c_1,w_1,c_2,w_2)$ where $c_1,c_2$ are odd vertex disjoint cycles and $w_1,w_2$
  are walks which combine a vertex $v_1$ of $c_1$ and a vertex $v_2$ of $c_2$. Then
  the result follows with the same way as in the corresponding
  proof of the \cite[Theorem 4.8]{BOO2}. \hfill$\square$

\section{Robustness and Generalized robustness on toric ideals}

In this section we study the robustness property
and the generalized robustness property in the general case
of toric ideals. In \cite{BOO2} Boocher et al proved that all robust toric ideals
which are generated by quadratics are graph ideals, see \cite[Corollary 5.3]{BOO2}.
As we can see in the next example, in the case of generalized robust toric ideals
this is not true.

\begin{ex1}{\rm
We consider the set $A=\{$(1,0,0,0,1), (0,1,1,1,1), (1,1,0,0,1), (0,0,1,1,1),
 (0,1,1,0,1), (1,0,0,1,1), (1,0,1,0,1), (0,1,0,1,1)$\}\subseteq \mathbb{N}^5.$
 We compute by CoCoA, that the corresponding
 toric ideal is $$I_A=<x_1x_2-x_3x_4,x_1x_2-x_5x_6,x_1x_2-x_7x_8>,$$ for more see \cite{COC}. The ideal $I_A$
 is quadratic and there is not a graph $G$ such that $I_A=I_G$, which means that $I_A$ is
 not a graph ideal. Otherwise, there exists a graph $G$ with three cycles of length four
 which have the edges $x_1$ and $x_2$ in common. Note that $x_1,x_2$ are vertex disjoint edges. This
 structure is impossible to happen for any simple graph $G$. Clearly for the toric ideal $I_A$
 we have that $Gr_A=<x_1x_2-x_3x_4,x_1x_2-x_5x_6,x_1x_2-x_7x_8,
 x_3x_4-x_5x_6,x_3x_4-x_7x_8,x_5x_6-x_7x_8>=\mathcal{M}_A$. By computations we
 check that $\mathcal{U}_A=\mathcal{M}_A$. It follows that the ideal $I_A$ is generalized robust.
}
\end{ex1}

As we saw in Corollary \ref{rob-gen} a robust graph ideal is also generalized robust.
Next we will see that this is a property not only for toric ideals of graphs, but for
 a random toric ideal as well. The main theorem of this section is that
robust ideals are generated by indispensable binomials.

In \cite[Theorem 2.12]{THO3} Charalambous et al,
described the indispensable elements of a toric ideal.

\begin{thm1}\cite{THO3}\label{indisp} The ideal $I_A$ is generated by indispensable binomials
if and only if the Betti $A-$degrees ${\bf{b}}_1,\ldots,{\bf{b}}_k$ are minimal
binomial $A-$degrees and $\beta_{0,{\bf{b}}_i}=1,\ \forall \ i=1,\ldots,k$.
\end{thm1}

Next, we remind some useful definitions, as they are presented in \cite{THO3},
in order to understand the above theorem.

Let $A\subset \mathbb{Z}^n$ be a vector configuration so that $\mathbb{N}A$ is pointed
and let $I_A$ be its corresponding toric ideal. A vector ${\bf{b}}\in \mathbb{N}A$ is
called a \textit{Betti A-degree} if $I_A$ has a minimal generating set containing an element
of $A-$degree ${\bf{b}}$. We define the \textit{A-graded Betti mumber} of $I_A$ as the
number of times that the vector ${\bf{b}}$ appears as the $A-$degree of a binomial in
a given minimal generating set of the ideal. From \cite{ST} we know that the Betti $A-$degrees
are independent of the choice of a minimal generating set of $I_A$.

Since the semigroup $\mathbb{N}A$ is pointed, we can partially order it with the relation:
$${\bf{c}}\geq{\bf{d}}\Longleftrightarrow \exists \ {\bf{e}}\in \mathbb{N}A: {\bf{c}}={\bf{d}}+{\bf{e}}.$$
Also for $I_A\neq\{0\}$, the minimal elements of the set
$\{\deg_A({\bf{x^u}}): {\bf{x^u}}-{\bf{x^v}}\in I_{A}\}\subset\mathbb{N}A$
with respect to $\geq$ are called \emph{minimal binomial $A-$degrees}.

For any ${\bf{b}}\in \mathbb{N}A$ the following ideal is defined:
$$I_{A,{\bf{b}}}=({\bf{x^u}}-{\bf{x^v}}\ : \ \deg_A({\bf{x^u}})=\deg_A({\bf{x^v}})\lneqq {\bf{b}})\subset I_A.$$

Next, we define the graph $G({\bf{b}})$. Based on this graph, an other graph ($S_{\bf{b}}$)
is defined, for more details see \cite{THO3}. The following construction plays a
key role in the proof of our main theorem.
\begin{def1}For a vector ${\bf{b}}\in \mathbb{N}A$ we define the graph $G({\bf{b}})$
to be the graph whose vertices are the elements of the fiber
$$\deg^{-1}_A({\bf{b}})=\{{\bf{x^u}}\ : \ \deg_A({\bf{x^u}})={\bf{b}}\}$$ and on the edge set
$$E(G({\bf{b}}))=\{\{{\bf{x^u}},{\bf{x^v}}\}\ : \ {\bf{x^u}}-{\bf{x^v}}\in I_{A,{\bf{b}}}\}.$$
\end{def1}
We consider the
complete graph $S_{\bf{b}}$, whose vertices are the connected components $G({\bf{b}})_i$ of $G({\bf{b}})$.
Let $T_{\bf{b}}$ be a spanning tree of $S_{\bf{b}}$. For every edge of $T_{\bf{b}}$ joining the
components $G({\bf{b}})_i$ and $G({\bf{b}})_j$ of $G({\bf{b}})$, we choose a binomial
${\bf{x^u}}-{\bf{x^v}}$ such that ${\bf{x^u}}\in G({\bf{b}})_i$ and ${\bf{x^v}}\in G({\bf{b}})_j$
correspondingly. Let $\mathcal{F}_{T_{\bf{b}}}$ the collection of these binomials.

\begin{prop1}\cite[Proposition 2.2]{THO3}\label{1} Let ${\bf{b}}\in \mathbb{N}A$. Every connected component
of $G({\bf{b}})$ is a complete subgraph. The graph $G({\bf{b}})$ is not connected if and
only if ${\bf{b}}$ is a Betti $A-$degree.
\end{prop1}

The following corollary, which results from the above definitions, is useful for the proof
of the main theorem of this article.

\begin{cor1}\label{difer} Let ${\bf{x^u}}-{\bf{x^v}}$ a minimal generator of $I_A$.
The vertices ${\bf{x^u}}$ and ${\bf{x^v}}$ of $G({\bf{b}})$ belong to different connecting components
of the graph $G({\bf{b}})$.
\end{cor1}
\textbf{Proof.} We assume that $$I_A=<B_1={\bf{x^{u_1}}}-{\bf{x^{v_1}}},\ldots,B_k={\bf{x^{u_k}}}-{\bf{x^{v_k}}}>,$$
where $\{B_1,\ldots,B_k\}$ is a minimal set of generators of $I_A$. Let ${\bf{b}}=\deg_A({\bf{x^{u_1}}}-{\bf{x^{v_1}}})$ be
the corresponding Betti $A-$degree.
We assume that the vertices ${\bf{x^{u_1}}}$ and ${\bf{x^{v_1}}}$ belong to the same connecting component
of the graph $G({\bf{b}})$. From Proposition \ref{1} we know that
every connected component of $G({\bf{b}})$ is a complete graph and therefore ${\bf{x^{u_1}}},{\bf{x^{v_1}}}$
is an edge of $G({\bf{b}})$. Thus, the binomial ${\bf{x^{u_1}}}-{\bf{x^{v_1}}}\in I_{A,{\bf{b}}}$, which means
that ${\bf{x^{u_1}}}-{\bf{x^{v_1}}}=\sum_ja_jM_jB_j$, where $j\in\{2,\ldots,k\}$. A
contradiction arises since the set $\{B_1,\ldots,B_k\}$ is minimal. \hfill$\square$

In \cite{THO3}
the authors proved the following theorems.

\begin{thm1}\cite[Theorem 2.6]{THO3}\label{Fminim} The set
$\mathcal{F}=\bigcup_{{\bf{b}}\in \mathbb{N}A}\mathcal{F}_{T_{\bf{b}}}$ is a minimal
generating set for the ideal $I_A$.
\end{thm1}

The converse of the above theorem is also true.

\begin{thm1}\cite[Theorem 2.7]{THO3}\label{minim-tree} Let $\mathcal{F}=\bigcup_{{\bf{b}}\in \mathbb{N}A}\mathcal{F}_{T_{\bf{b}}}$
be a minimal generating set of the ideal $I_A$.
The binomials of $\mathcal{F}_{T_{\bf{b}}}$ determine a spanning tree
$T_{\bf{b}}$ of $S_{\bf{b}}$.
\end{thm1}

The next proposition will be used in the sequel.

\begin{prop1}\cite[Proposition 2.4]{THO3}\label{mindeg} An $A-$degree ${\bf{b}}$ is minimal binomial
$A-$degree if and only if every connected component of $G({\bf{b}})$ is a singleton.
\end{prop1}

Before presenting the main theorem, we put forward the following lemma, which was proved
in \cite{BOO2} in the special case of toric ideals of graphs. We
remark that there is no difference in the general case of a toric
ideal. By $\mu(I_A)$ we denote the number of minimal generators of the ideal.

\begin{lem1}\label{Boo} Let $I_A$ be a robust toric ideal. Then there is no term of an
element of $\mathcal{U}_A$ which divides a term of another element of $\mathcal{U}_A$.
\end{lem1}
\textbf{Proof.} Let $I_A$ be a robust toric ideal.
By definition, the set $\mathcal{U}_A$ is a minimal generating set for the ideal.
Since the affine semigroup $\mathbb{N}A$ is pointed,
the graded Nakayama Lemma applies that all minimal system of generators of $I_A$
have the same cardinality. Thus $|\mathcal{U}_A|=\mu(I_A)$.
Also we remark that from the definition of robustness, it follows that
$\mu(\textrm{in}_\prec I_A)=\mu(I_A),\ \textrm{for all term orders} \prec $
and therefore $$|\mathcal{U}_A|=\mu(\textrm{in}_\prec I_A),\ \textrm{for all term orders} \prec .$$
We will prove the contrapositive. We suppose that $\mathcal{U}_A$ contains at
least two binomials $f_1=m_1-m_2$ and $f_2=n_1-n_2$ such that the term
$m_1$ divides the term $n_1$. Since the toric ideal $I_A$ is prime,
there exists a variable $x$ which divides the monomial $m_1$ but not the monomial $m_2$.
Taking $\prec$ to be the lex term order with $x$ first, it follows that
$(\textrm{in}_\prec f_1)|(\textrm{in}_\prec f_2)$. Therefore
$$\mu(\textrm{in}_\prec I_A)< |\mathcal{U}_A|,$$ a contradiction arises.
\hfill$\square$

We are ready now to continue with the proof of
the main theorem of this section. In the following theorem
we are proving that the robust toric ideals have a unique minimal system of generators.

\begin{thm1}\label{main} Let $I_A$ be a robust toric ideal. Then $I_A$ is generated by indispensable binomials.
\end{thm1}
\textbf{Proof.} Let $I_A$ be a robust toric ideal. From definition the ideal $I_A$
is minimally generated by its universal Gr\"obner basis. Let
$$M=\{B_1={\bf{x^{u_1}}}-{\bf{x^{v_1}}},B_2={\bf{x^{u_2}}}-{\bf{x^{v_2}}},\ldots,B_k={\bf{x^{u_k}}}-{\bf{x^{v_k}}}\}$$
be the universal Gr\"obner basis of $I_A$, which is also a
minimal generating set for the ideal.
We will prove that the binomials $B_i$ are indispensable generators of $I_A$
for every $i=1,\ldots,k$. From Theorem \ref{indisp}, it is enough to prove
that the Betti $A-$degrees ${\bf{b}}_1,\ldots,{\bf{b}}_k$ are minimal
binomial $A-$degrees and $\beta_{0,{\bf{b}}_i}=1,\ \forall \ i=1,\ldots,k$.
As we saw before, we know that the Betti $A-$degrees
are independent of the choice of the minimal generating set of $I_A$ and
therefore if we choose at random one of them, it appears as a degree of an
element of the set $M$ as well.

Firstly, we prove that the Betti $A-$degrees are minimal.
We consider the minimal generator $B_1={\bf{x^{u_1}}}-{\bf{x^{v_1}}}$
and let ${\bf{b}}_1=\deg_A({\bf{x^{u_1}}}-{\bf{x^{v_1}}})$ be its Betti $A-$degree. Let $G({\bf{b}}_1)$
be the corresponding graph. From
Proposition \ref{mindeg}, the Betti $A-$degree ${\bf{b}}_1$ is minimal if and only
if every connected component of $G({\bf{b}}_1)$ is a singleton. Let $G({\bf{b}}_1)_i$ be
a connected component of $G({\bf{b}}_1)$ which is not a singleton. Without loss
of generality we can assume that the vertex ${\bf{x^{u_1}}}$ belongs to $G({\bf{b}}_1)_i$
and let ${\bf{x^{w}}}$ be an other vertex of $G({\bf{b}}_1)_i$. Note that from Corollary \ref{difer},
the vertex ${\bf{x^{v_1}}}$ does not belong in $G({\bf{b}}_1)_i$ and let $G({\bf{b}}_1)_j$ be
its connected component. We consider the corresponding tree $T_{{\bf{b}}_1}$. Then
from Theorem \ref{minim-tree}, the edge $e_1=\{{\bf{x^{u_1}}},{\bf{x^{v_1}}}\}$ of the
graph $S_{{\bf{b}}_1}$
is an edge of $T_{{\bf{b}}_1}$. We replace this edge
by the edge $e=\{{\bf{x^w}},{\bf{x^{v_1}}}\}$ in the tree $T_{{\bf{b}}_1}$.
Obviously, we have a new spanning tree of the graph $S_{{\bf{b}}_1}$ and from
Theorem \ref{Fminim} the set
$$M'=\{B_2={\bf{x^{u_2}}}-{\bf{x^{v_2}}},\ldots,B_k={\bf{x^{u_k}}}-{\bf{x^{v_k}}},B_{k+1}={\bf{x^{w}}}-{\bf{x^{v_1}}}\}$$
is a minimal generating set of $I_A$. We remark that there is no element of $M'$
which contains the term ${\bf{x^{u_1}}}$. Otherwise the monomial ${\bf{x^{u_1}}}$ appears
at least twice in two different elements of the set $M$, which is the universal Gr\"obner bases
of $I_A$, a contradiction arises from Lemma \ref{Boo}. We consider the binomial
$B_1={\bf{x^{u_1}}}-{\bf{x^{v_1}}}$ which belongs to $I_A$. Since
$M'$ is a minimal generating set, the $B_1$ can be written as
a linear combination of the elements of $M'$. We have that
$${\bf{x^{u_1}}}-{\bf{x^{v_1}}}=\sum_{j=2}^{k+1}a_jm_jB_j,\ \textrm{where}\
a_j\in \mathbb{K}\ \textrm{and}\ m_j\ \textrm{monomials of}\ \mathbb{K}[{\bf{x}}].$$
As a result, there exists $j\in\{2,\ldots,k,k+1\}$ such that ${\bf{x^{u_j}}}| {\bf{x^{u_1}}}$.
Note that $j\neq k+1$, otherwise ${\bf{x^{w}}}| {\bf{x^{u_1}}}$
which is impossible, since ${\bf{x^{w}}}\neq {\bf{x^{u_1}}}$ and
$\deg_A({\bf{x^{w}}})=\deg_A({\bf{x^{u_1}}})$.
Therefore, there is a term of an element of $M=\mathcal{U}_A$, which divides
a term of another element of $\mathcal{U}_A$. From Lemma \ref{Boo}, we have a
contradiction. So, we conclude that
every connected component of $G({\bf{b}}_1)$ is a singleton and therefore
the Betti $A-$degree ${\bf{b}}_1$ is minimal.

It remains to prove that $\beta_{0,{\bf{b}}_i}=1,\ \forall \ i=1,\ldots,k$.
We assume that there exists $i\in\{1,\ldots,k\}$ such that $\beta_{0,{\bf{b}}_i}\geq2$.
In other words, there are at least two elements of $M$ with degree ${\bf{b}}_i$
and let them be $B_n={\bf{x^{u_n}}}-{\bf{x^{v_n}}}$ and $B_m={\bf{x^{u_m}}}-{\bf{x^{v_m}}}$.
As we proved before, every connecting component of
$G({\bf{b}}_i)$ is a singleton, thus the graph $G({\bf{b}}_i)$ has at least four
connecting components; $\{{\bf{x^{u_n}}}\},\{{\bf{x^{v_n}}}\},\{{\bf{x^{u_m}}}\}$ and $\{{\bf{x^{v_m}}}\}$.
Note that none of the above connecting components coincides with each other, otherwise
we have a contradiction from the Lemma \ref{Boo}.
If we look at the corresponding tree $T_{{\bf{b}}_i}$, two of its edges are
$e_1=\{{\bf{x^{u_n}}},{\bf{x^{v_n}}}\}$ and $e_2=\{{\bf{x^{u_m}}},{\bf{x^{v_m}}}\}$.
Since a tree is a connected graph, then there exists a path which joins the edges $e_1,e_2$.
Therefore, at least one of the vertices $\{{\bf{x^{u_n}}}\},\{{\bf{x^{v_n}}}\}$
appears in an other one edge of $T_{{\bf{b}}_i}$, different from $e_1$. This means
that at least one of the monomials ${\bf{x^{u_n}}}$ or ${\bf{x^{v_n}}}$ appears as
a monomial term of
an other minimal generator of $M=\mathcal{U}_A$. A contradiction arises from Lemma \ref{Boo}.
Thus $\beta_{0,{\bf{b}}_i}=1,\ \forall \ i=1,\ldots,k$ and the theorem follows.
\hfill$\square$

The converse of the above theorem is not true, as we can see in the following remark.

\begin{rem1}{\rm
The indispensability of the minimal generators of a toric ideal,
is not a necessary condition for an ideal to be robust.
 For this claim we consider the following graph.

\begin{figure}[here]
\begin{center}
\psfrag{A}{$e_1$}\psfrag{B}{$e_2$}\psfrag{C}{$e_3$}\psfrag{D}{$e_4$}\psfrag{E}{$e_5$}
\psfrag{F}{$e_6$}\psfrag{G}{$e_7$}\psfrag{H}{$e_8$}
\psfrag{I}{$e_9$}\psfrag{J}{$e_{10}$}\psfrag{K}{$e_{11}$}\psfrag{L}{$e_{12}$}
\includegraphics{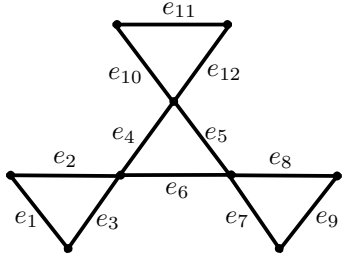}
\caption{An indispensable ideal which is not robust}
\label{Figure 4}
\end{center}
\end{figure}

The corresponding toric ideal $I_G$ has six minimal generators. These are:
$$B_1=e_2e_3e_5-e_1e_4e_6, B_2=e_4e_5e_{11}-e_6e_{10}e_{12},B_3=e_4e_7e_8-e_5e_6e_9,$$
$$B_4=e_2e_3e_{10}e_{12}-e_1e^2_4e_{11},B_5=e_2e_3e_7e_8-e_1e^2_6e_9,B_6=e_7e_8e_{10}e_{12}-e_9e^2_5e_{11}.$$

From Corollary \ref{uniqueminm} all the above minimal generators of $I_G$ are indispensable.
We consider the walk $w=(e_2,e_1,e_3,e_6,e_7,e_9,e_8,e_5,e_{12},e_{11},e_{10},e_4)$ of the graph $G$.
From Theorem \ref{primitive-graph} the corresponding binomial
$B_w=e_2e_3e_7e_8e_{10}e_{12}-e_1e_4e_5e_6e_9e_{11}$
belongs to the Graver basis of $I_G$. We note that there are two sinks of $w$ in distance one and
therefore $w$ is not strongly primitive. From Theorem \ref{minimal} the binomial $B_w$ is not minimal.
Theorem \ref{M=Gr} implies that the ideal $I_G$ is not generalized robust. Therefore it is not robust.
}
\end{rem1}

We note that although the robust ideals are generated by indispensable binomials, this
does not happen in the case of generalized robust ideals, see Example \ref{ex1}.

By Theorem \ref{main} it follows the next corollary, in which we see that the property
of robustness for a toric ideal,
implies the generalized robustness property for it.
\begin{cor1} \label{eis}Let $I_A$ be a toric ideal. If $I_A$ is robust then it is generalized robust.
\end{cor1}

Obviously, a necessary condition for the converse statement of the above corollary,
is the uniqueness of the minimal system of generators of the toric ideal.

\begin{cor1} Let $I_A$ be a toric ideal, such that it has a unique
minimal system of generators. The ideal $I_A$ is robust if and only if $I_A$ is
generalized robust.
\end{cor1}

\begin{rem1}\label{remm}{\rm We note that the equality
$M_A= Gr_A$ for robust toric ideals still remains an open problem, as it has
been mentioned in \cite{BOO2}.
We remark that the intersection of all minimal Markov bases $M_A$ of a toric ideal, is called the
indispensable subset of the universal Markov basis $\mathcal{M}_A$ and is denoted by $S_A$.
Obviously, the following inclusions hold: $$S_A\subseteq \mathcal{M}_A\subseteq Gr_A.$$ In
\cite{THO2} the authors gave a complete algebraic characterization for the
elements of the sets $S_A$ and $\mathcal{M}_A$. As it follows from the previous theorem,
in the case of robust toric ideals we have that $S_A=\mathcal{M}_A$. An equivalent interesting
question for robust toric ideals or generalized robust toric ideals
is the equality $\mathcal{M}_A=Gr_A$. The interest for this problem is enhanced by the fact that
toric ideals for which the universal Gr\"obner basis coincides with their Graver basis
have important properties, as for example the equality between the Gr\"obner complexity
and the Graver complexity of the ideal, see \cite{HM}. Also, in the case of robust toric ideals
the following inclusions hold:
$$\mathcal{C}_A\subseteq (S_A=\mathcal{M}_A=M_A=)\ \mathcal{U}_A\subseteq Gr_A.$$
There are examples of families of ideals whose set of circuits is equal with the Graver basis. For example, Sturmfels
proved this property for toric ideals defined by unimodular matrices, see \cite[Proposition 8.11]{ST}
and Villareal proved it for those defined by balanced matrices, see \cite{VI2}. We know that
toric ideals of graphs which are complete intersection are also circuit ideals, which means that
every minimal generator of the ideal is a circuit, see \cite[Theorem 5.1]{TT3}. Thus, for those
toric ideals we have that $\mathcal{C}_A=\mathcal{M}_A$. By Theorem \ref{M=Gr},
we get an other family of ideals, i.e. robust toric ideals of graphs which are complete intersection,
for which $\mathcal{C}_A=Gr_A$.
}
\end{rem1}

\end{document}